\documentclass[letterpaper, 10 pt, conference]{ieeeconf} 
\IEEEoverridecommandlockouts
\usepackage{cite}
\usepackage{amsmath,amssymb,amsfonts}
\usepackage{algorithmic}
\usepackage{graphicx}
\usepackage{textcomp}
\usepackage{xcolor}

\usepackage[left=1.9cm, right=1.9cm, top=2.54cm]{geometry}
\title{\LARGE \bf
	The Multi-objective Dynamic Traveling Salesman Problem: Last Mile Delivery with Unmanned Aerial Vehicles Assistance 
}

\author{Ben Remer, Andreas A. Malikopoulos, \textit{Senior Member, IEEE}
\thanks{This research was supported by the Delaware Energy Institute (DEI).}
	\thanks{The authors are with the Department of Mechanical Engineering, University of Delaware, Newark, DE 19716 USA (email:  \texttt{bremer@udel.edu}; \texttt{andreas@udel.edu}.) }%
}

\begin{document}
	
	\bibliographystyle{unsrt}
	
	\maketitle
	\thispagestyle{empty}
	\pagestyle{empty}

	
	
	\indent
	\begin{abstract}
		In this paper, we present an approach to optimizing the last-mile delivery route of a truck using coordination with unmanned aerial vehicles (UAVs). First, a traveling salesman problem is formulated to determine the truck's route. Then, a scheduling problem is formulated to determined the routes for the UAVs. A genetic algorithm is used to solve these problems, and simulated results are presented.
	\end{abstract}
	
	\section{Introduction}
In a rapidly urbanizing world, we need to make fundamental transformations
in how we use and access transportation. We are currently witnessing an increasing integration of our energy, transportation, and cyber networks, which, coupled with the human interactions, is giving rise to a new level of complexity in transportation \cite{Malikopoulos2015}. 
As we move to increasingly complex transportation systems \cite{Malikopoulos2016c},
new control approaches \cite{Malikopoulos2015b, Malikopoulos2018d} are needed to optimize the impact on system
behavior of the interaction between vehicles at different applications. 

	Unmanned aerial vehicles (UAVs) are becoming increasingly available to be applied in civilian life, from drone racing to filming of events. In the past few years, applying UAVs to perform better delivery services has become a research topic, partially spurred by Amazon's ``Air Amazon" delivery concept announced in 2013. Potentially, a delivery truck with a team of UAVs could increase both the time and energy efficiency of a last mile delivery service (that is, the portion of the delivery going from the distribution center to the final destination)\cite{Ponza2016}.
	
	Last mile delivery is a research focus, considered the least efficient part of the delivery. There are many different approaches in the literature, such as crowd sourcing of deliveries \cite{Wang2016} in which the goal is to increase efficiency by decreasing missed deliveries. Another related approach in the literature are \textit{Green Vehicle Routing Problems} \cite{Lin2014}\cite{Yao2015}. These problems may sometimes try to increase efficiency of last mile delivery by reducing energy usage. Use of UAVs coordinate with a ground vehicle to perform last-mile delivery has attracted some attention lately \cite{Murray2015},\cite{Dorling2017}. Typically, the UAV is taken to be a quadrotor, which is also what we consider in this paper. The majority of research efforts have focused on optimizing the UAV itself \cite{Dorling2017}, including some relevant efforts into the control schemes of UAVs in conjunction with changing physical parameters due to the package delivery \cite{Quadrotors2018a}, than optimizing the last mile delivery as a whole.
	
	A slight version of this problem was introduced in 2015 \cite{Murray2015},\cite{Mathew2015}. In their approach, the authors make the assumption that the UAV and the truck will only meet up and depart at some truck delivery node -- albeit Murray and Chu \cite{Murray2015} permit the UAV to leave and enter the depot separately from the truck, Mathew, Smith, and Waslander \cite{Mathew2015} do not. This assumption was later relaxed in \cite{Bari2017a}, permitting the UAV and the truck to rendezvous along any edge the truck is traveling. However, there is still the assumption that the truck is parked during this action. This limitation is unnecessary, as landing a quadrotor on a moving target has been explored with some promising results \cite{Xu2016}. In this paper, we relax this assumption to explore the potential benefits. To the best of our knowledge this problem has not been addressed in the literature to date. 

	In several research efforts reported in literature, generally a variant of the traveling salesman problem (TSP) has been formulated. In \cite{Mathew2015} the problem was addressed as a \textit{heterogeneous delivery problem}, while in \cite{Murray2015} the problem was addressed as the \textit{flying sidekick traveling salesman problem}. In the TSP, there is a collection of cities, each of which must be visited exactly once by a salesman, and a depot that the salesman starts and ends at. The goal of the problem is to find the optimal route between these cities. These classes of problems are NP-Hard. In the literature, it is shown that our problem is unique from the standard TSP, as we have multiple agents. In addition, in some formulations, and this is the case in this paper, we have delivery points (or cities) that can only be visited by a delivery truck since the packages might be too heavy for a UAV, and some that can be visited by either the delivery truck or the UAV\cite{Ponza2016},\cite{Murray2015}. In addition, we formulate the problem such that the UAV is energy and range constrained. Additionally, we expand upon the literature such by formulating the problem by allowing \textit{n} number of UAVs.
	
	

	The rest of this paper is organized as follows. In Section II, we cover the problem formulation of the truck route. In Section III, we address the plausible rendezvous and departure points for the UAV. In Section IV, we present the mathematical formulation of the scheduling problem associated with the UAVs the algorithms for solving the problem. In Section V, we present our simulation results, and finally, in Section VI we draw concluding remarks and discuss potential avenues for future research.

	\section{Mathematical Formulation}
	\begin{figure}
		\includegraphics[width=.45\textwidth]{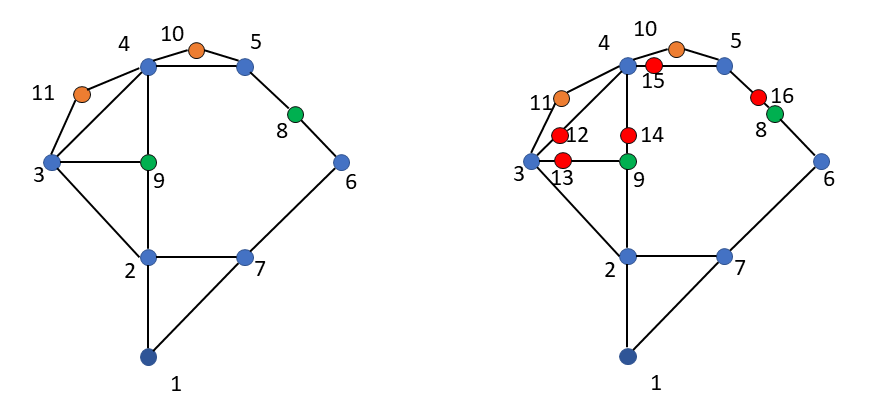}
		\caption{ Left is the original graph $\mathcal{G}$. Right is the graph modified into $\mathcal{G}'$ by adding rendezvous nodes}
	\end{figure}

	Before discussing the problem at hand, we briefly review some basic notions of graph theory. We use $\mathcal{G(N,A)}$ to describe a graph with a set of arcs (or edges) $\mathcal{A}$ connecting a set $\mathcal{N}$ of $N$ nodes (or vertices), i.e., $\mathcal{N}=\{1,\dots,N\},N\in \mathbb{N}$. Each node is indexed by $i \in \mathcal{N}$. We take a system of a truck and $n$ UAVs cooperating to make deliveries. The network of the truck-UAVs system is represented by the graph $\mathcal{G(N,A)}$. 
	
	An example of the network is shown in Fig. 1. We reserve $i=1$ to be the depot node that the truck and UAVs must start and end together. In the example illustrated in Fig. 1, the set depot node 1 is displayed in dark blue color. We define four disjoint sets $\mathcal{D},\mathcal{H},\mathcal{S},$ and $\mathcal{R}$, in which each individual node belongs to.  The set $\mathcal{D} \subset \mathcal{N}$ is a set of delivery nodes that can be serviced by the UAVs, although the truck may still service them. In the example illustrated in Fig. 1, the set $\mathcal{D}$ has two nodes displayed in orange color. The set $\mathcal{H} \subset \mathcal{N}$ is the set of delivery nodes that must be reached by the truck only. The set $\mathcal{H}$ has two nodes displayed in green color. The set $\mathcal{S} \subset \mathcal{N}$ is the set that includes the topology nodes (e.g., street intersections). The set $\mathcal{S}$ has six nodes displayed in blue color. Finally, the set $\mathcal{R} \subset \mathcal{N}$ is the set of possible rendezvous and deploy nodes for the UAVs. Initially, this set is empty but an algorithm can identify and locate these nodes.
	Fig. 1 shows the transformation from the original graph (left) to the graph (right) where the set $\mathcal{R}$ of rendezvous and deploy nodes is included. A second algorithm can detect possible $sortie$. A sortie is defined as a 3-node route $i-j-k$ that a UAV can take, where $i, k\in \mathcal{R}$ and $j\in\mathcal{D}$. The set of all sorties that can visit a node $i$ is denoted by $\mathcal{F}_i$, while an individual sortie is denoted $s_i$. 
	
	Finally, we find the shortest path between each member of the set $\mathcal{H}\cup\mathcal{D}\cup\mathcal{R}\cup{S}$, and use this to create our graph $\mathcal{G''(N'',A'')}$. 
	As the dynamic element in our problem is the change in mass, the relative lengths between nodes do not change, we can safely make this formulation. This will allows us to constrain our solver to visit each node no more than once, but not prevent us from revisiting street nodes.
	
	In our problem formulation, $r$ denotes our decisions variable, which is an array of nodes that describes the route of the truck. We also denote with $x_{ij}$ the edge that is selected.  For instance, if $x_{12}=1$, then the truck travels from node 1 to node 2. This, of course, occurs when node 1 is followed by node 2 in $r$. We define $r_i$ to be the $i$th member of $r$ and $r_{-1}$ to be the last member. Thus, we have
	\begin{equation}
	x_{ij}=
	\begin{cases}
	1 \text{ if $r_k=j$ and $r_{k-1}=i$ for some $r_k\in r$}\\
	0 \text{ otherwise}.
	\end{cases}
	\end{equation}
	
	We are interested in both time and energy of the entire truck-UAVs system, and thus we formulate a multi-objective optimization problem as follows
	
	\begin{equation}
	\begin{aligned}
	\min \big(\alpha E(r)+(1-\alpha)T(r)\big), ~
	0\leq\alpha\leq 1,
	\end{aligned}
	\end{equation}
	s.t.
	
	\begin{equation}
	\label{TruckDelivery}
	\sum_{i \neq j}x_{ij} =1~~ \forall j \in H,
	\end{equation}
	\begin{equation}
	\label{Start}
	r_i=1,
	\end{equation}
	\begin{equation}
	\label{End}
	r_{-1}=1,
	\end{equation}
	\begin{equation}
	\label{Flow}
	\sum_{i}x_{ij} \leq 1 ~~\forall j,
	\end{equation}
	\begin{equation}
	\label{Sortie}
	\sum_{i \in \mathcal{F}_j}\delta(s_i) \geq 0 ~~\forall j \in {\mathcal{D} \notin r},
	\end{equation}
	where,
	\begin{equation}
	\label{delta}
	\delta(s_i)=
	\begin{cases}
	1 \text{ if rendezvous nodes of $s_i \in r$},\\
	0 \text{ otherwise}.
	\end{cases}
	\end{equation}
	In our problem formulation above, $\alpha$ is a weighting parameter, $E(r)$ is the energy costs and $T(r)$ is the total time to complete deliveries. A breakdown of the costs $E(r)$ and $T(r)$ is presented in Section II. The delivery constraint (\ref{TruckDelivery}) ensures truck deliveries are made, while the constraints (\ref{Start}) and (\ref{End}) ensure the truck starts and ends at the depot node. The constraint (\ref{Flow}) implies that we can only visit a node once; however, recall this constraint applies on the transformed graph, so it is possible to revisit the same street nodes, if needed, due to the connectivity of the original graph. Finally, the constraint (\ref{Sortie}) states that for each delivery being made by the UAVs and not the truck, there is at least one possible sortie. This permits possible solutions for the scheduling portion of the problem. It does by making use of a function $\delta(s)$ which probes if a sortie's end points are within the route.
	
	$Edge$ $Cost:$ For edge costs, we are concerned with time and energy. First, we consider the case of the truck and the the UAV in flight. Note that we define the edge costs generically as to handle the case of the UAVs with the truck, and with the changing mass of both the UAV and truck as packages get delivered. Additionally, we need to store what the route on the original graph makes up the route that this edges transits.
	
	$Node$ $Label:$ To store all necessary information for a given node, we create a label for it. To label each node, we use a tuple $(i,Q,M_p,(x,y,z))$. $i$ is a unique index, $Q$ is the set the node belongs to (either $\mathcal{H}$, $\mathcal{D}$, $\mathcal{S}$, or $\mathcal{R}$), $M_p$ is the mass of the package that must be delivered there, if there is no package to be delivered $M_p = 0$, and finally $(x,y,z)$ are the euclidean coordinates of the delivery where $z$ is the height of the delivery point, or the height of the truck.
	
	In our modeling framework, we make the following assumptions.
	
	\textit{Assumption 1}: The UAV can fly in a direct route path to its destination node.
	
	\textit{Assumption 2}: There is no wind velocity and changes of elevations are negligible.
	
	\textit{Assumption 3}: The effect the package being delivered has on the quadrotor's frontal area and coefficient of drag is negligible.
	
	\textit{Assumption 4}: The UAVs travel from the truck, to a delivery point, and back to the truck. The UAVs land on top of, and take off from the top of the truck. Additionally, it is assumed that immediately after landing on the truck the UAVs have their batteries swapped and packages loaded by an automated system.
	
	\textit{Assumption 5}: Outside of initial and final acceleration and braking, the truck moves at the maximum permitted and constant speed along each edge.

	Assumption 1 ensures a plausible path for the UAV to take reliably. Physically, this means that any buildings or landmarks are shorter than the height the UAV is flying at, or at least not in the way for whatever reason. Assumption 2 simplifies the energy model for the UAV and the truck, and implies that the costs to transverse an edge are the same in both directions as long as the mass is the same. Assumption 3 is an approximation as the package's contribution is small to the frontal area and coefficient of drag. Additionally, this is a practical assumption as it is extremely unlikely for these effects to be measured in practice. Assumption 4 is simply an assumption of execution, which is necessary for an answer to be achieved (for instance, if the UAV were permitted to ride in the truck, the truck's frontal area would not change with the UAV being on board. If the UAV would be to exposed partially, the frontal area would change some. This is a decision that needs to be made in design). Assumption 5 implies that the truck moves at a predictable and deterministic way while along an edge. This is an idealized behavior of a ``bang-bang" controller which will minimize the time to cross an edge, disregarding any real world factors that could disrupt the truck's velocity, e.g. traffic.
	
	\subsection{Edge Energy Costs of the Truck in Transit}

	The truck's energy cost can be broken into the energy associate with the acceleration and cruising at a velocity\cite{Rios-Torres2,Rios-Torres2015,Kamal2013,Kamal2011}. The mass of the truck will change as deliveries are made, monotonically decreasing throughout, for the TSP portion of our problem as we will ignore the mass fluctuation from the UAV riding the truck. The cost the truck endures from the UAV riding the truck will be taken into account during the scheduling portion of the problem
	
	\begin{equation}
	e_{u,ij}=\int_{accel}\hat{a}_{ij}(c_0+c_1v+c_2v^2)dt,
	\end{equation}
	where
	\begin{equation}
	\hat{a}_{ij}=-(1)/(2M_{i})\cdot C_{d}\rho_aAv^2-\mu g+u_h,
	\end{equation}
	subject to:
	$u_{brake}\leq u_h \leq u_{acc},$
	$u_{brake} < 0 < u_{acc},$ and
	$0\leq v \leq V_{ij},$
	where $e_{u,ij}$ is the energy the truck consumes for a given input acceleration $u$, which for simplicity we will take to be the maximum or minimum value. $\hat{a}$ is the equivalent acceleration, $u_{brake}$ and $u_{acc}$ are vehicle parameters, $A$ is the frontal area of the truck, $M_{ij}$ is the mass of the truck across edge $ij$, $C_{d}$ is the coefficient of drag, $\rho_a$ is the density of air,  $\mu$ is a friction coefficients, and $u_h$ is the acceleration input. This model has the energy usage cut off when decelerating from a high velocity. We need also to include the energy associated with cruising the truck, which is give as follows
	\begin{equation}
	e_{cruise,ij}=\int_{cruising}b_0+b_1V_{ij}+b_2V_{ij}^2+b_3V_{ij}^3dt
	\end{equation}
	where $b_0$, $b_1$, $b_2$, and $b_3$ are vehicle parameters. Therefore, the total energy the truck expends going from node $i$ to node $j$ is
	\begin{equation}
	e_{t,total,ij}=e_{u,ij}+e_{cruise,ij}.
	\end{equation}
	
	This leads to our energy cost problem formulation
	\begin{equation}
	\begin{aligned}
	E(r)={} & w_1\sum_{i \in N, j \in N, i\neq j} e_{t,total,ij}(x^1_{ij})\\
	&+w_2\sum_{k=m-c+1}^m\sum_{i \in N, j \in N, i\neq j} e^k_{d,total,ij}(x^k_{ij})
	\end{aligned}
	\end{equation}
	where $w_1$ and $w_2$ are weights comparing the relative value of the energy spent of the UAV and the energy spent of the truck. For instance, this could correspond to the cost of gasoline for the truck and electricity for the UAV. Our energy model for the UAV gives power consumption in Watt$\cdot$sec, but we convert to kWh before proceeding. The energy of the truck is measured in milliliters of fuel. Thus, a reasonable value for $w_1$ could be around 0.000747 (based of of 2.87 dollars/gallon of fuel) and a reasonable value for $w_2$ would be around 0.12 (12 cents/kwh). Note that while $w_2 > w_1$, the truck energy term is expected to be significantly larger than the UAV energy term for a normal range of physical values.
	
	\subsection{Time Costs}
	
	The time costs to transverse an edge can be computed based on the velocity, accelerations, and the distance. In doing so, we apply Assumption 5 that outside of initial acceleration and final deceleration, the truck moves at a constant velocity equal to the maximum velocity permitted along that edge
	\begin{equation}
	\begin{aligned}
	t_{t,ij}={} & (d_{ij}-v_{ij}^2/(2u_{acc})+v_{ij}^2/(2u_{brake}))/v_{ij} \\
	& +v_{ij}/u_{acc}-v_{ij}/u_{brake} \text{ if $i \in \mathcal{H}, j \in \mathcal{H}$}
	\end{aligned}
	\end{equation}
	\begin{equation}
	\begin{aligned}
	t_{t,ij}={} & (d_{ij}+v_{ij}^2/(2u_{brake}))/v_{ij} \\
	& -v_{ij}/u_{brake} \text{ if $i \notin \mathcal{H}, j \in \mathcal{H}$}
	\end{aligned}
	\end{equation}
	\begin{equation}
	\begin{aligned}
	t_{t,ij}={} & (d_{ij}-v_{ij}^2/(2u_{acc})))/v_{ij} \\
	& +v_{ij}/u_{acc} \text{ if $i \notin \mathcal{H}, j \in \mathcal{H}$}
	\end{aligned}
	\end{equation}
	\begin{equation}
	\begin{aligned}
	t_{t,ij}={} & d_{ij}/v_{ij} \text{ if $i \notin \mathcal{H}, j \notin \mathcal{H}$}.
	\end{aligned}
	\end{equation}
	
	Hence,
	
	\begin{equation}
	T(r)=\sum_{j \in N}\sum_{i \in N, i\neq j}t_{t,ij}(x_{ij}).
	\end{equation}
	
	Note that the time for the UAV to cross any edge is not in our final equation. This is because the truck and UAV are constrained to rendezvous and must start and end together.
	
	\subsection{Edge Energy Costs of the UAV in Flight}
	From Assumption 1, we can generate a path for the UAV $k$ between nodes. We can compute the energy cost of traveling from a node $i$ to a node $j$, denoted $e^k_{d,total,ij}$, as a sum of the energy cost associated with ascending to that height, denoted $e^k_{a,ij}$, descending from that height, denoted $e^k_{d,ij}$, flying the distance between the two nodes, denoted $e^k_{t,ij}$, and the costs associated with performing the rendezvous maneuver (if applicable).
	We use the energy model described in \cite{Liu2017} in which the authors reported results within 10\% of tested quadrotors. We stress that our future analysis is not dependent on this model, but requires an energy model that is dependent on the changing mass of the quadrotor. Furthermore, we note that there are some particular quadrotor concepts reported in the literature \cite{Driessens2013} that this model might not be appropriate for as they adopt unusual rotor configurations that are not compatible with assumptions made in the energy model. Finally, we simplify this model further with Assumption 2. The energy cost of ascending with this model is

	\begin{equation}
	\label{AscentE}
	\begin{aligned}
	e^k_{a,ij}={}&(z-z_i)/V^k_a[d^k_2(m^k_{ij}g)^{3/2}+k^k_1m^k_{ij}g\cdot\\
	&\cdot[V^k_{a}/2+\sqrt{(V^k_{a}/2)^2+m^k_{ij}g/(k^k_2)^2}]],
	\end{aligned}
	\end{equation}
	where $k^k_1$ and $k^k_2$ are physical parameters of the UAV $k$ that can be found experimentally, $z_i$ is the height of the node $i$, $g$ is acceleration due to gravity, $m^k_{ij}$ is the mass of the UAV during the ascent (note: this value changes with the mass of the package being carried), and $V^k_a$ is the velocity during ascent. The descent is formulated identically in (3), however $V^k_d$ is the descent velocity and is negative
	
	\begin{equation}
	\label{DescentE}
	\begin{aligned}
	e^k_{d,ij}={}&(z_j-z)/V^k_d[d^k_2(m^k_{ij}g)^{3/2}+k^k_1m^k_{ij}g\cdot\\
	&\cdot[V^k_{d}/2+\sqrt{(V^k_{d}/2)^2+m^k_{ij}g/(k^k_2)^2}]].
	\end{aligned}
	\end{equation}
	
	Finally, we look at the energy costs associated with the UAV traveling between two nodes. We consider the power loss from drag, the power to hover, and the profile power
	
	\begin{equation}
	\label{TransverseE}
	\begin{aligned}
	e^k_{t,ij}=\int_{transverse}(P^k_{hover,ij}+P^k_{par,ij}+P^k_{p,ij})dt\\
	P^k_{hover,ij}=d^k_1(T^k_{ij})^{3/2} \\
	P^k_{par,ij}=d^k_4V_{ij}^k3 \\
	P^k_{p,ij}=d^k_2(T^k_{ij})^{3/2}+d^k_3(V^k_{ij}\cos(\alpha))^2(T^k_{ij})^{1/2} \\
	T^k_{ij}=\sqrt{(m^k_{ij}g-d^k_5((V^k_{ij}\cos(\alpha))^2)^2+(d^k_4V_{ij}^{k2})^2},
	\end{aligned}
	\end{equation}
	where $d^k_1$, $d^k_3$,$d^k_4$, and $d^k_5$ are physical parameters of the UAV $l$. $\alpha$ is the angle of attack. Changes in frontal area and drag coefficients with the change in the package being carried are negligible, if Assumption 3 holds.
	
	Taking these factors into consideration, the total energy cost for the UAV to transverse an edge $ij$ while flying can be written as
	\begin{equation}
	e^k_{d,total,ij}=e^k_{t,ij}+e^k_{a,ij}+e^k_{d,ij} \text{ if } i \in \mathcal{D} \lor j \in \mathcal{D}.
	\end{equation}
	
	As expected, this is a function of the UAV's physical parameters, the time to transverse the edge, and the package weight. We have to additional cases, covering the situation of the UAV riding along the truck. Recall that this is only for UAV in flight, in other words one of the ends nodes is a delivery node. For when the UAV is on the truck, we model it as an additional cost of energy with the same parameters of the truck but with the same mass of the UAV. 
	
	\subsection{Updating Mass}
	As deliveries are made, the mass of the truck and the UAV change. The mass of the truck goes down monotonically as deliveries are made-recall that even while the UAV is docked on the truck we handle the energy spent to move it separately from the truck. The mass of the UAV fluctuates with each delivery.
	
	We denote $M_0^*$ the starting mass of the truck-UAV-packages system and $M_0^f$ the final mass of the system at the depot node. For every other node $i$, we denote $M_i$ to be the mass at that node, after any deliveries, departures, and rendezvouses are made. Finally, we define $m_i$ to be the mass of the delivery at node $i$. If no delivery is to be made at node $i$, $m_i=0$.
	
	We can define the mass at each node as
	
	\begin{equation}
	M_j=\sum_{i \in N} x_{ij}M_i-m_j.
	\end{equation}
	
	We separately consider the mass of the UAV $k$ across edge $ij$ to be
	
	\begin{equation}
	m^k_{ij}=
	\begin{cases}
	M_{UAV,k}+m_j \text{ if $j \in \mathcal{R}$}\\
	M_{UAV,k} \text{ otherwise}.
	\end{cases}
	\end{equation}
	
	\section{Rendezvous and UAV Departure Points}
	The first step we will take to solve this problem, is add plausible rendezvous and departure points to the original graph. From our battery constraint (\ref{Battery}) we can deduce that any acceptable rendezvous and departure point $j \in \mathcal{R}$ to be paired with a delivery node $i \in \mathcal{D}$ for UAV $k$ has a battery constraint
	$e^k_{batt} > e^k_{d,total,ij}.$
	Taking (\ref{AscentE})-(\ref{TransverseE}) and the battery constraint, the maximum distance the UAV can cover under this constraint can be solved for. Due to Assumptions 1 and 2, this is the same distance in any direction, as there will be no obstacles to avoid and no wind. We call this $R_{max}$. Finally, for each delivery node $D_i$ we just need to find every edge within this radius. For each edge, we create a node along that edge (the precise location of the node is not defined, just that it is along that edge) and add an edge between that pair of nodes and the accessible delivery nodes.  After finding the rendezvous nodes, we can then find each set of sorties.

For each delivery node, we simply check each possible pair of edges between that node and a rendezvous node (note that in this case a ``pair" can contain the same edge twice). If the sum of the energy costs of those two edge is less than the battery requirement of the UAV, it is added to set of sorties.
	
	
	\section{The Unmanned Aerial Vehicle Scheduling Subproblem}
	Before proceeding, we review some scheduling notation and verbiage and discuss what it means in the context of our problem. In scheduling, there are machines, in our case UAVs, that perform jobs, in our case deliveries. We say our machines are parallel machines, as they perform jobs independently of each other. Finally, we are interested in finding a series of times $(t_1,t_2,\dots,t_3)$ for each machine in which it perform each job $(j_1,j_2,\dots,j_N)$ that was assigned to it in the scheduling process. Each time for job $i$ on machine $j$ is constrained by a release time $r_{ij}(t_{i-1},geometry,V_0,V_f)$, which is the earliest time that the machine can start the job. It is further constrained by a deadline $d_{ij}(geometry,V_0,V_f)$ and a processing time $p_{ij}(t_i,geometry,V_0,V_f)$. While in some scheduling problems $preemption$ is allowed, it is not permitted here. 
	
	After our algorithm has defined a route for the truck, a scheduling problem must be done for the UAVs. The optimal solution is one that minimizes the penalty function
	\begin{equation}
	\label{ScheduleCost}
	J=\sum_{k=1}^n J^k,
	\end{equation}
	where,
	$J^k(T^k_1,\dots,T^k_n)=\sum_i^n C^k(T_i),$
	where $C^k(T_i)$ is the cost of UAV $k$ doing job $i$ at time $T$. It is derived from the geometry of the problem, the physical parameters of the truck and UAVs, and the energy models previously discussed. Additionally, a dynamical model for the UAVs must be adopted to fully describe the route for the energy model. We adopt a simple double integrator for the sake of simplicity. While more sophisticated models for UAVs exist and methods such as minimum snap exist for determining paths, that remains outside the scope of this paper. $(T_1,\dots,T_n)$ describe the time in which a UAV $k$ starts the deliveries $1,\dots,n$. We will then sum the penalty functions for each UAV, and add it to the previous solution. Finally, we will rerun the scheduling subproblem for any previously found solutions that are now better than the previously identified optimal solution.
	
	\subsection{Solving the Scheduling Problem}
	We start by discretizing the problem, and for each sortie to perform the action on UAV $k$ a we calculate the total cost and end time $T^k_f$ for a given start time $T^k_s$. If the cost is greater than the battery constraint (\ref{Battery}), then that time is excluded. We also ignore sorties that are not plausible due to having an end node not in our truck solution $r$. We apply the constraints
	$T^k_{i,f}<T^k_{i+1}.$
	In other words, a UAV cannot start a job until it has finished the previous one. In addition, we have the previously mentioned battery constraint
	\begin{equation}
	\begin{aligned}
	\label{Battery}
	\sum_{i=1}^N e_{d,total,ij}\cdot x^k_{ij} + \sum_{i=1}^N e_{d,total,ji}\cdot x^k_{ij} \leq e^k_{batt}\\
	\forall j \in D \text{ and } \forall k \mid 2\leq k \leq n+1.
	\end{aligned}
	\end{equation}
	We then seek a solution that minimizes the scheduling problem subject to the time constraint for a given set of $N$ jobs on UAV $k$ $(j_1,j_2,\dots,j_N)$.
	
	Two popular algorithms for solving similar problems are Simulated Annealing (SA) and Genetic Algorithms (GA). SA in particular has been used to solve problems relating to a truck-UAV tandem making deliveries \cite{Mathew2015}. In this paper, we will take the opportunity to explore a genetic algorithm.

	\section{Simulation Results}
	The algorithms discussed in this paper were written in MATLAB, and tested on a road network that resembles a scaled up version of the University of Delaware's Scaled Smart City Lab (UDSSC) \cite{Malikopoulos2018b}.
	All delivery points were randomly generated.
	A total of 3 randomly generated UAVs were used. 4 truck deliveries (location and mass) were randomly generated and 8 UAV deliveries (location and mass) were randomly generated. These results were then compared to the results of having no UAV assistance.
	
	Average results from simulations are shown in Table I, with $\alpha$ skewed aggressively towards time at 0.9. Assisted results refer to those achieved with UAV assistance, and  unassisted refers to results achieved without UAV assistance (for the same randomly generated delivery points). A 20.77\% improvement is achieved on average, representing significant savings that can be achieved.
	
	\begin{table}
		\label{results}
		\caption{Average Values of Cost Function}
		\begin{center}
			\renewcommand{\arraystretch}{1.4}
			\begin{tabular*}{.405\textwidth}{|l|c|c|r|}
				\hline			
				Assisted  & Unassisted & Improvement & \% Improvement \\
				\hline
				372.94 & 473.02 & 100.07 & 20.77\% \\
				\hline
			\end{tabular*}
		\end{center}
	\end{table}
	
	\section{Concluding Remarks}
	In this paper, we presented an approach to optimizing the last-mile delivery route of a truck using coordination with UAVs. First, a traveling salesman problem was formulated to determine the truck's route. Then, a scheduling problem was formulated to determined the routes for the UAVs. A genetic algorithm is used to solve these problems, and simulated results are presented.
		Providing assistance to a delivery truck, while adding complexity, can improve the efficiency of last mile delivery. We provided a problem formulation that is flexible enough to handle multiple and unique UAVs assisting the truck, while providing a trade off between time and energy when present. We adopted a unique two-problem algorithm to branch and bound the solution and provide fast and good results.
		

	\bibliographystyle{IEEEtran}
	\bibliography{references}

\begin{thebibliography}{10}

\bibitem{Malikopoulos2015}
A.~A. Malikopoulos.
\newblock Centralized stochastic optimal control of complex systems.
\newblock In {\em Proceedings of the 2015 European Control Conference}, pages
  721--726, 2015.

\bibitem{Malikopoulos2016c}
A.~A Malikopoulos.
\newblock A duality framework for stochastic optimal control of complex
  systems.
\newblock {\em IEEE Transactions on Automatic Control}, 61(10):2756--2765,
  2016.

\bibitem{Malikopoulos2015b}
A.~A. Malikopoulos, V.~Maroulas, and J.~XIong.
\newblock A multiobjective optimization framework for stochastic control of
  complex systems.
\newblock In {\em Proceedings of the 2015 American Control Conference}, pages
  4263--4268, 2015.

\bibitem{Malikopoulos2018d}
A.~A. Malikopoulos, C.D. Charalambous, and I.~Tzortzis.
\newblock The average cost of markov chains subject to total variation distance
  uncertainty.
\newblock In {\em Systems \& Control Letters}, volume 120, pages 29--35, 2018.

\bibitem{Ponza2016}
A.~Ponza.
\newblock {Optimization of Drone-Assisted Parcel Delivery}.
\newblock page~80, 2016.

\bibitem{Wang2016}
Y.~Wang, D.~Zhang, Q.~Liu, F.~Shen, and L.~Hay.
\newblock {Towards enhancing the last-mile delivery : An effective
  crowd-tasking model with scalable solutions}.
\newblock {\em Transportation Research Part E}, 93:279--293, 2016.

\bibitem{Lin2014}
Canhong Lin, K.L. Choy, G.T.S. Ho, S.H. Chung, and H.Y. Lam.
\newblock {Survey of Green Vehicle Routing Problem : Past and future trends}.
\newblock {\em Expert Systems with Applications}, 41:1118--1138, 2014.

\bibitem{Yao2015}
Enjian Yao, Zhifeng Lang, Yang Yang, and Yongsheng Zhang.
\newblock {Vehicle routing problem solution considering minimising fuel
  consumption}.
\newblock {\em IET Intelligent Transport Systems}, 9(5):523--529, 2015.

\bibitem{Murray2015}
Chase~C. Murray and Amanda~G. Chu.
\newblock {The flying sidekick traveling salesman problem: Optimization of
  drone-assisted parcel delivery}.
\newblock {\em Transportation Research Part C: Emerging Technologies},
  54:86--109, 2015.

\bibitem{Dorling2017}
Kevin Dorling, Jordan Heinrichs, Geoffrey~G. Messier, and Sebastian
  Magierowski.
\newblock {Vehicle Routing Problems for Drone Delivery}.
\newblock {\em IEEE Transactions on Systems, Man, and Cybernetics: Systems},
  47(1):70--85, 2017.

\bibitem{Quadrotors2018a}
Animesh~Kumar Shastry, Mahathi~T Bhargavapuri, Mangal Kothari, and
  Soumya~Ranjan Sahoo.
\newblock {Quaternion Based Adaptive Control for Package Delivery using
  Variable-Pitch Quadrotors}.
\newblock In {\em Indian Control Conference}, 2018.

\bibitem{Mathew2015}
Neil Mathew, Stephen~L. Smith, and Steven~L. Waslander.
\newblock {Planning Paths for Package Delivery in Heterogeneous Multirobot
  Teams}.
\newblock {\em IEEE Transactions on Automation Science and Engineering},
  12(4):1298--1308, 2015.

\bibitem{Bari2017a}
P~D~I Bari, L~Caggiani, M~Marinelli, M~D Orco, and M~Ottomanelli.
\newblock {En-route truck-drone parcel delivery for optimal vehicle routing
  strategies}.
\newblock 12:253--261, 2017.

\bibitem{Xu2016}
L~Xu and H~Luo.
\newblock {Towards Autonomous Tracking and Landing on Moving Target *}.
\newblock 2016.

\bibitem{Rios-Torres2}
J.~Rios-Torres and A.~A Malikopoulos.
\newblock {Automated and Cooperative Vehicle Merging at Highway On-Ramps}.
\newblock {\em IEEE Transactions on Intelligent Transportation Systems},
  18(4):780--789, 2017.

\bibitem{Rios-Torres2015}
J.~Rios-Torres, A.~Malikopoulos, and P.~Pisu.
\newblock {Online Optimal Control of Connected Vehicles for Efficient Traffic
  Flow at Merging Roads}.
\newblock {\em IEEE Conference on Intelligent Transportation Systems,
  Proceedings, ITSC}, 2015-Octob:2432--2437, 2015.

\bibitem{Kamal2013}
A.~S. Kamal, M.~Mukai, J.~Murata, and T.~Kawabe.
\newblock {Model Predictive Control of Vehicles on Urban Roads for Improved
  Fuel Economy}.
\newblock 21(3):831--841, 2013.

\bibitem{Kamal2011}
M.~A.S. Kamal, Masakazu Mukai, Junichi Murata, and Taketoshi Kawabe.
\newblock {Ecological vehicle control on roads with up-down slopes}.
\newblock {\em IEEE Transactions on Intelligent Transportation Systems},
  12(3):783--794, 2011.

\bibitem{Liu2017}
Zhilong Liu, Raja Sengupta, and Alex Kurzhanskiy.
\newblock {A Power Consumption Model for Multi-rotor Small Unmanned Aircraft
  Systems}.
\newblock In {\em International Conference on Unmanned Aircraft Systems
  (ICUAS)}, pages 310--315, 2017.

\bibitem{Driessens2013}
S~Driessens and P~E~I Pounds.
\newblock {Towards a more efficient quadrotor configuration}.
\newblock {\em 2013 IEEE/RSJ International Conference on Intelligent Robots and
  Systems}, pages 1386--1392, 2013.

\bibitem{Malikopoulos2018b}
A~Stager, L.~Bhan, A.~A Malikopoulos, and L.~Zhao.
\newblock A scaled smart city for experimental validation of connected and
  automated vehicles.
\newblock In {\em 15th IFAC Symposium on Control in Transportation Systems},
  pages 130--135, 2018.

\end{thebibliography}
	
\end{document}